\documentclass[12pt,a4paper,twoside,final,notitlepage, leqno]{article}
\usepackage[english]{babel}
\usepackage[T1]{fontenc}  
\usepackage{epsfig, graphicx, amssymb}
\usepackage{amsmath,amsthm,epsfig,amsfonts,bbm}  
\usepackage{float}
\usepackage{color}
\newcommand{\monitem}{ \smallskip \noindent $\bullet$ \quad  } 
\newcommand{\moneq}{\vspace*{-7pt} \begin{equation} \displaystyle } 
\newcommand{\moneqstar}{\vspace*{-6pt} \begin{equation*} \displaystyle } 
\newcommand{\monendstar}{\vspace*{-6pt} \end{equation*}   }
\newcommand{\monend}{\vspace*{-7pt} \end{equation}   }

\newcommand{\R}{\mathbb{R}}

\newcommand{\Z}{\mathbb{Z}}

\def\section*#1{}

\usepackage{fancyhdr}
\renewcommand{\headrulewidth}{0pt}

\begin{document} 

\fancypagestyle{plain}{ \fancyfoot{} \renewcommand{\footrulewidth}{0pt}}
\fancypagestyle{plain}{ \fancyhead{} \renewcommand{\headrulewidth}{0pt}}

~

  \vskip 2.1 cm

\centerline {\bf \LARGE  Curious convergence properties }

 \bigskip 

\centerline {\bf \LARGE  of lattice Boltzmann schemes }

 \bigskip 

\centerline {\bf \LARGE   for diffusion with acoustic scaling }

 \bigskip  \bigskip \bigskip

\centerline { \large    Bruce M. Boghosian$^{a}$, Fran\c{c}ois Dubois$^{bc}$, Benjamin Graille$^{b}$,}

\smallskip  \centerline { \large   Pierre Lallemand$^{d}$ and Mohamed-Mahdi Tekitek$^{e}$}

\smallskip  \bigskip 

\centerline { \it  \small   
$^a$  Dpt. of Mathematics, Tufts University, Bromfield-Pearson Hall, Medford, MA 02155, U.S.A.} 

\centerline { \it  \small   
$^b$   Dpt. of Mathematics, University Paris-Sud,  B\^at. 425, F-91405  Orsay, France.} 

\centerline { \it  \small   
$^c$    Conservatoire National des Arts et M\'etiers, LMSSC laboratory,  F-75003 Paris, France.} 

\centerline { \it  \small  $^d$   Beijing Computational Science Research Center, 
Haidian District, Beijing 100094,  China.}

\centerline { \it  \small  $^e$   Dpt. Mathematics, Faculty of Sciences of Tunis, 
University Tunis El Manar, Tunis, Tunisia. } 


\bigskip  

\centerline { 24  March 2017 
{\footnote {\rm  \small $\,$ Contribution published in {\it Communications in Computational Physics}, 
 volume 23, issue~4, pages 1263-1278, April 2018, doi: 10.4208/cicp.OA-2016-0257.
Edition 06 March 2018.}}}

 \bigskip 
 {\bf Keywords}: artificial compressibility method, Taylor expansion method.  

 {\bf PACS numbers}:  
02.70.Ns, 
05.20.Dd, 
47.10.+g 

\bigskip  
\noindent {\bf \large Abstract} 

\noindent 
We consider  the D1Q3 lattice Boltzmann scheme with an acoustic scale
for the simulation of diffusive  processes. 
When the mesh is refined while holding the diffusivity constant, we first obtain asymptotic convergence.
When the mesh size tends to zero, however, this convergence breaks down in a curious fashion, and
we observe qualitative discrepancies from analytical solutions of the heat equation. 
In this work, a  new asymptotic analysis is derived to explain this phenomenon using the 
Taylor expansion method, and a partial differential equation of acoustic type 
is obtained in the asymptotic limit.
We show that the error between the D1Q3 numerical solution and 
a finite-difference approximation of this acoustic-type partial differential equation
tends to zero in the asymptotic limit.
In addition, a wave vector analysis of this asymptotic regime
demonstrates that the dispersion equation has nontrivial complex eigenvalues, 
a sign of underlying propagation phenomena,
and a portent of the unusual convergence properties mentioned above.

\bigskip \bigskip   \newpage \noindent {\bf \large    1) \quad  Introduction }   

\fancyhead[EC]{\sc{ B. Boghosian, F. Dubois,  B. Graille, P. Lallemand,  M.M. Tekitek }} 
\fancyhead[OC]{\sc{Curious convergence  of lattice Boltzmann schemes for diffusion... }} 
\fancyfoot[C]{\oldstylenums{\thepage}}

\noindent 

Lattice Boltzmann models are simplifications of the continuum Boltzmann 
equation obtained by discretizing in both physical space and velocity space.  
The discrete velocities $v_i$ retained typically correspond to lattice vectors 
of the discrete spatial lattice.  That is, each lattice vertex $x$ is linked 
to a finite number of neighboring vertices by lattice vectors $v_i \, \Delta t$.
A particle distribution $ \, f \, $ is therefore parametrized by its components in each of the discrete velocities, the vertex $x$ of the spatial lattice, and the discrete time $t$.  A time step of a classical lattice Boltzmann scheme~\cite{LL00} then contains two steps:

\noindent
{\it (i)} a relaxation step where distribution $ f $ at each vertex $x$ 
is locally modified into a new  distribution $ f^* $, and

\noindent
 {\it (ii)} an advection step based  on the method of characteristics 
as an  exact time-integration operator.  We employ the  multiple-relaxation-time 
approach introduced by d'Humi\`eres~\cite{DDH92}, wherein 
the local mapping $ \, f \longmapsto  f^*  \, $ is described by a nonlinear 
diagonal  operator in a space of moments, as detailed in Section~2.

\noindent 
In~\cite{DL09}, we have 
studied the asymptotic expansion of various 
lattice Boltzmann schemes with multiple-relaxation times 
for different applications.  We used the so-called acoustic scaling, in which
the ratio  $ \, \lambda \equiv \Delta x / \Delta t \,$ is kept fixed. 
In this manner, we demonstrated
the possibility of approximating diffusion processes described by the heat equation.

\noindent 
In his very complete work, Dellacherie~\cite{De14} has described unexpected results in simulations 
for advection-diffusion processes.  In this contribution, 
we endeavor to explain those results by studying the 
convergence of the D1Q3 lattice Boltzmann scheme when we try to approximate 
a pure diffusion process. 

\noindent 
We begin this paper by recalling some fundamental algorithmic aspects of the D1Q3 lattice Boltzmann
scheme in Section~2.  Then, in Section~3 we describe a first illustrative numerical experiment.
In Section~4 we present a new convergence analysis, followed by another numerical experiment in
Section~5, in which the D1Q3 lattice Boltzmann scheme is studied far from the usual values of
its parameters.  Finally, a wave vector analysis is proposed in Section~6.

\bigskip \bigskip   \noindent {\bf \large    2) \quad  
Diffusive D1Q3 lattice Boltzmann scheme }   

\noindent 
In this work, we consider the so-called D1Q3 lattice Boltzmann scheme in one spatial dimension. 
The  spatial step $ \, \Delta x > 0 \, $ is given, and each node  $ \, x \,$
is an  integer multiple of this spatial step~:  $ \, x \in  \Z  \, \Delta x $.
The time step  $ \, \Delta t > 0 \,$ is likewise given, and each discrete time $ \, t \,$
is an  integer multiple of  $ \, \Delta t$.
We adopt so-called acoustic scaling (see {\it e.~g.}, \cite{UAS10}), so the numerical velocity
associated with the mesh,
\moneq \label{lambda} 
 \lambda \equiv {{\Delta x}\over{\Delta t}} ,
\monend
is a constant independent  of the spatial step $ \, \Delta x $. 
A  particle distribution 
\moneqstar  
f \, \equiv \,  \big( f_+(x,\, t)  \,,\,  f_0(x,\, t) \,,\,   f_- (x,\, t) \big) 
\monendstar 
is given at the initial step $ \, t = 0 $.  Its value at subsequent times
is determined by the multiple-relaxation-time version~\cite{DDH92} 
of the lattice Boltzmann equation.

\monitem 
Moments are introduced at each step of space and time according to the relations 
\moneq \label{moments}  
\rho = f_+ + f_0 + f_-  \,,  \quad 
J  = \lambda \, ( f_+ - f_- )   \,,  \quad  
e  = \lambda^2  \, ( f_+ - 2 \, f_0 + f_- ) \, . 
\monend
These may be thought of as the densities of mass, momentum, and an energy-like quantity, respectively.
Eq.~(\ref{moments}) can be recast in matrix form as follows:
\moneqstar 
m  \,   \equiv \,  \begin{pmatrix} \rho \cr J  \cr e  \end{pmatrix} 
\,=\, M \, f  \,   \equiv \,  M \,  \begin{pmatrix} f_+ \cr f_0 \cr f_-   \end{pmatrix} \,, 
\monendstar 
where $ \, M \,$ is the invertible matrix,
\moneq \label{matrice-M}  
M  \,   =  \,  \begin{pmatrix} 1  & 1 & 1 \cr 
\lambda & 0  & - \lambda \cr \lambda^2  & -2 \, \lambda^2  &  \lambda^2  \end{pmatrix}  \, . 
\monend
%

\monitem 
The equilibrium values of the moments are defined by the relations:
\moneq \label{equilibre}  
\rho^{\rm eq}  = \rho  \,,  \quad  J^{\rm eq}  = 0  \,,  \quad  
e^{\rm eq}  =  \alpha \, {{\lambda^2}\over{2}} \, \rho 
\monend
where $\, \alpha \,$ a non-dimensional constant. Then the relaxation step transforms 
the pre-collision moments $ \, m \,$ into new post-collision  moments $ \, m^* $ as follows: 
\moneq \label{moments-collision}  
\rho^*  = \rho  \,,  \quad  
J^{*}  = J + s_J \, ( J^{\rm eq} - J )   \,,  \quad  
 e^{*}  = e + s_e  \, ( e^{\rm eq} - e )  \, ,
\monend
where $s_J$ and $s_e$ are relaxation parameters.  There is no analogous parameter for $\rho$ because
the collisions are constrained to conserve mass.  In our numerical experiments, we have chosen $ \, s_e = 1.5 $, and below we shall explain in some detail how we tuned the relaxation parameter $ \, s_J \, $ for the momentum $ \, J $.

\monitem 
The time iteration of the scheme is defined in terms of the particle distribution. 
We first transform  the post-collision moments $ \, m^* \, $ into a
post-collision particle distribution: 
\moneqstar   
f^* \,=\, M^{-1} \, m^* \, . 
\monendstar 
Second, we iterate the algorithm forward in time.  The particle distribution is conserved 
along the characteristic directions of velocities $ \, v_+ = \lambda $,
$ \, v_0 = 0 \, $ and  $ \, v_- = -\lambda \,$ respectively:
\moneq \label{schema-temps}  
\left\{  \begin{array} {l}
\displaystyle f_+ (x, \, t+ \Delta t) \,=\, f_+^* (x-\Delta x , \, t )  \cr  
f_0 (x, \, t+ \Delta t) \,\,=\, f_0^* (x , \, t ) \cr   
f_- (x, \, t+ \Delta t) \,=\, f_-^* (x+\Delta x , \, t )  \, . 
\end{array} \right.  \monend
%

\monitem 
In~\cite{DL09}, we have analyzed several lattice Boltzmann models 
with the Taylor expansion method, including the present one defined by
Eqs.~(\ref{moments},\ref{matrice-M},\ref{equilibre},\ref{moments-collision},\ref{schema-temps}).
The hypothesis used was that the numerical velocity $ \, \lambda \,$ defined in Eq.~(\ref{lambda}),
and the relaxation coefficients $ \, s_J \,$ and $ \, s_e \,$ remain constant 
as the spatial step $ \, \Delta x \, $ tends to zero. 
Then the  conserved variable~$ \, \rho \, $  satisfies (at least formally!)
a diffusion partial differential equation: 
\moneq \label{edp-diff} 
{{\partial \rho}\over{\partial t}}   \, - \, \mu \,\,  
{{\partial^2 \rho}\over{\partial x^2}}  \,=\, {\rm O}(\Delta x^2)  \, ,
\monend
where the diffusion coefficient  $ \mu  \, $ is given by the relation
\moneq \label{mu} 
\mu  \,\equiv\, {{4+\alpha}\over{6}} \, \sigma  \, \lambda \,  \Delta x \,,\quad  
\sigma \, \equiv \, \Big({{1}\over{s_J}} - {{1}\over{2}} \Big) \, . 
\monend
The coefficient $ \, \sigma \,$ is known as the ``H\'enon parameter''
in reference to the pioneering work of 
H\'enon~\cite{He87}. This lattice Boltzmann scheme 
is demonstrably stable under the condition: 
\moneqstar 
-4 \, < \, \alpha \, < \, 2 \, . 
\monendstar

\bigskip    \noindent {\bf \large    3) \quad  
A first numerical experiment   }   

\begin{figure}  [H]     \centering
 \includegraphics [width=.60 \textwidth, angle=-90] {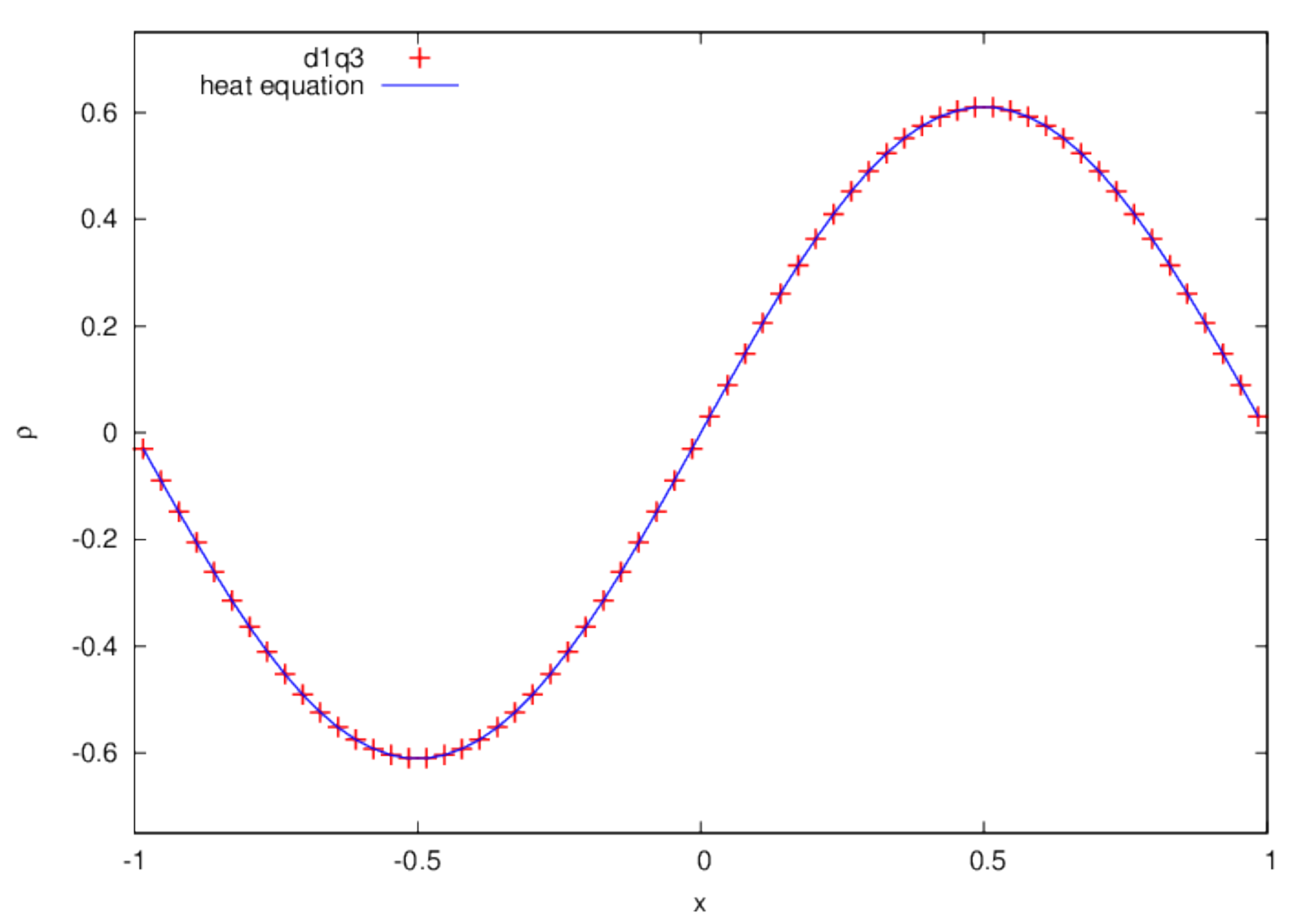} 
\caption{Sine wave  at time $ \, T = 5 \, $
computed with 64 mesh cells. }
  \label{sine-time5} \end{figure}

\begin{figure}  [H]     \centering
 \includegraphics [width=.60 \textwidth, angle=-90] {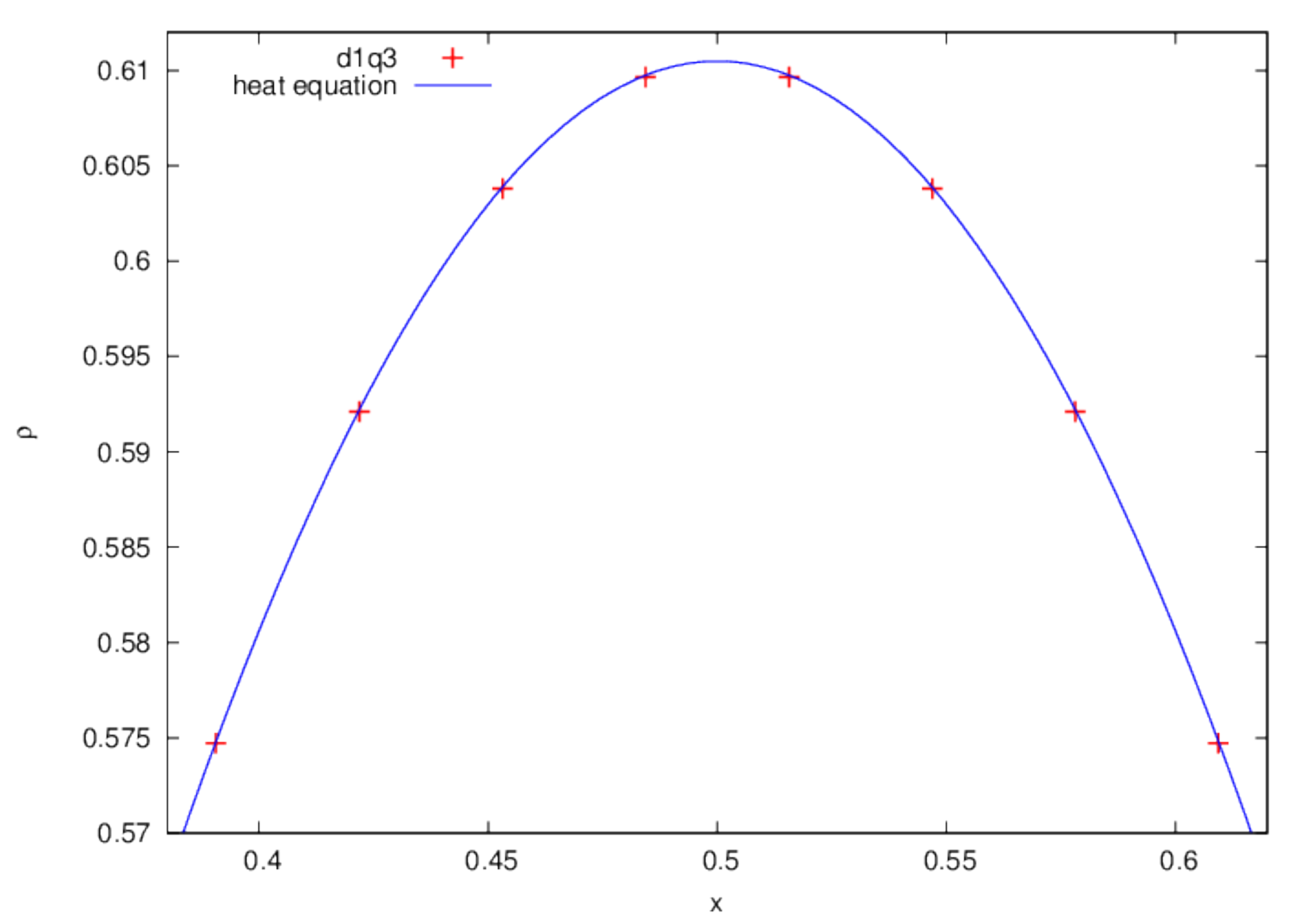} 
\caption{Sine wave  computed with 64 mesh cells at time $ \, T = 5 \, $; detail.}
  \label{sine-time5-detail2} \end{figure}

\begin{figure}  [H]     \centering
 \includegraphics [width=.60 \textwidth, angle=-90] {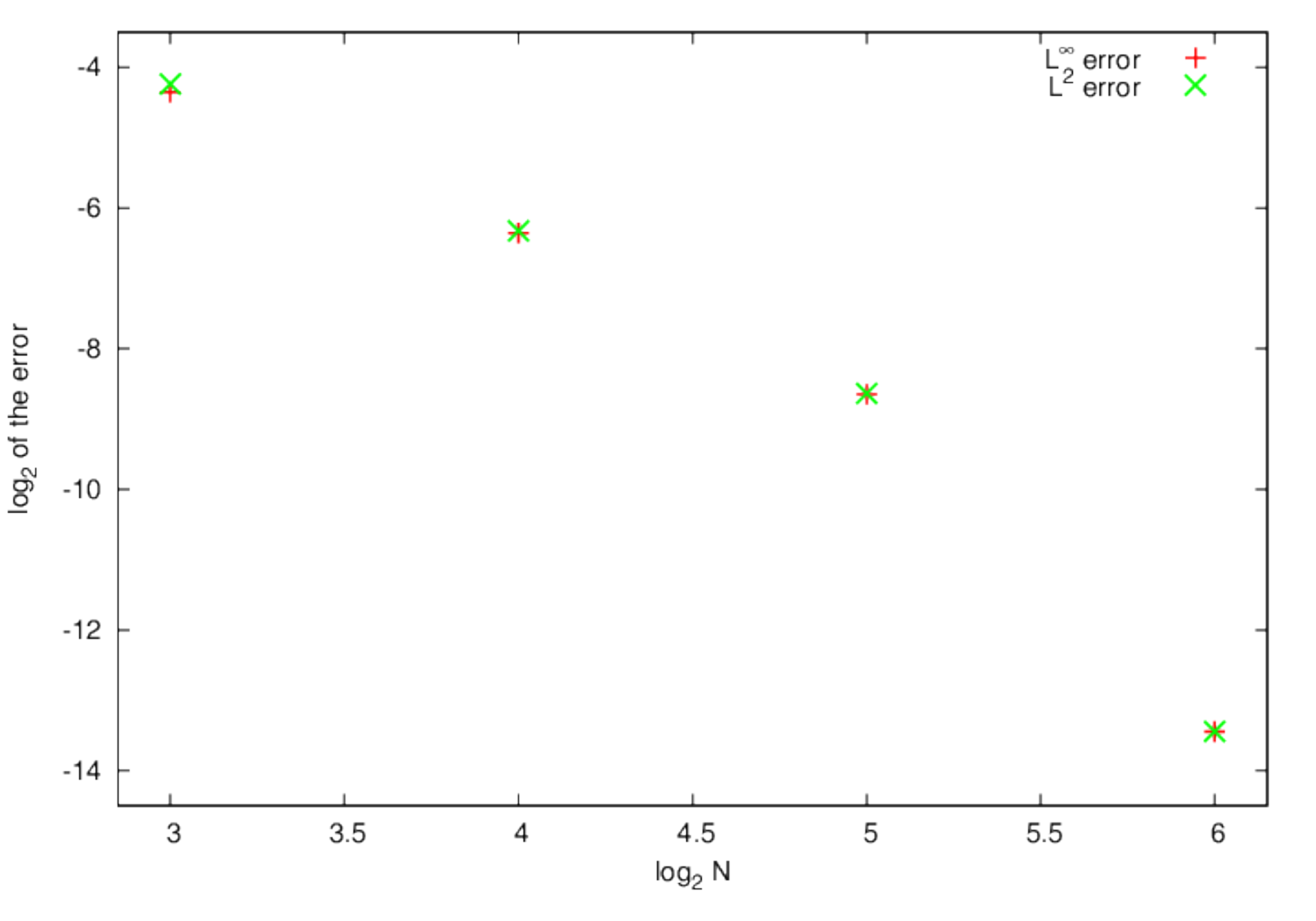} 
\caption{Numerical convergence to the diffusion equation~(\ref{edp-diff}).}
  \label{sine-time5-converge} \end{figure}

\noindent
In this section, we consider an elementary analytic test case, namely the diffusion of a sine 
wave.  We suppose that the initial condition for Eq.~(\ref{edp-diff})
satisfies 
\moneq \label{sinus-initial} 
 \rho_0(x) =  \sin (\pi \ x)  \,, \quad   -1 \leq x \leq 1  \, .  
\monend
The other moments $\, J \,$ and $ \, e \,$ are taken at equilibrium at $ \, t = 0 $. 
With periodic boundary conditions, the exact solution of Eqs.~(\ref{edp-diff},\ref{sinus-initial}) is 
\moneqstar 
 \rho(x,\, t) =  \sin ( \pi \, x ) \, \exp(-\mu \, \pi^2 \, t)  \, . 
\monendstar 
We performed several numerical computations with the following choice of parameters:
$ \, \lambda  =  1 $, $ \, \alpha  =  1 \, $ and $ \, \mu = 0.01 $. 
%
The spatial step varied from $ \, \Delta x = {1\over4} \,$ up to $ \, \Delta x = {1\over32} $. 
The results for a final time $ \, T = 5 \,$ 
are presented in Figs.~\ref{sine-time5} through \ref{sine-time5-converge}.

\noindent
It should be noted that the results are remarkably converged even for these relatively coarse meshes. 
For the most refined mesh used (64 mesh points, $ \Delta x = {1\over32} $), 
the numerical results are almost indistinguishable from the exact solution, 
as presented in two successive magnifications in Fig.~\ref{sine-time5-detail2}. 

\bigskip \bigskip   \noindent {\bf \large    4) \quad  
An alternative convergence analysis }   

\noindent 
We now imagine that we wish to approximate the diffusion equation, Eq.~(\ref{edp-diff}),
using the D1Q3 lattice Boltzmann model described previously.  
We suppose that the diffusion  coefficient  $ \, \mu  \, $ 
is fixed and that the mesh size $ \, \Delta x \,$
tends to zero.  Then from Eq.~(\ref{mu}), the relaxation parameter $ \,s_J \,$ 
can no longer be fixed and tends to zero according to the asymptotic prescription
\moneq \label{s-mu} 
s_J \,=\, {{4+\alpha}\over{6 \, \mu }} \,  \lambda \,  \Delta x \,+\, {\rm O}(\Delta x^2)  \, . 
\monend
The hypothesis used for  deriving the diffusion model, Eqs.~(\ref{edp-diff},\ref{mu}) 
is now violated, because
the relaxation parameter $ \,s_J \,$ is no longer a  constant.  Rather, it follows the asymptotic form
\moneq \label{s-dx} 
s_J \,=\, s_0 + s_1 \, \Delta x +  {\rm O} ( \Delta x^2 )   +  \cdots ,
\monend
as suggested by one of us in earlier work on lattice-gas automata~\cite{BB}. 
Moreover, we have 
\moneq \label{developpement-s} 
s_0 \, = \, 0 \,,\quad s_1 =  {{4+\alpha}\over{6 \, \mu }} \, 
\monend
for the case described in Eq.~(\ref{s-mu}).  Then the differential equation obtained in the asymptotic limit is no longer the diffusion equation Eq.~(\ref{edp-diff}), as discussed in the following proposition. 

\bigskip  \monitem {\bf Proposition 1. \quad An asymptotically acoustic model }

\noindent 
We consider the D1Q3 lattice Boltzmann scheme defined by
Eqs.~(\ref{moments},\ref{matrice-M},\ref{equilibre},\ref{moments-collision},\ref{schema-temps}).
We make the hypothesis that the numerical velocity 
$ \, \lambda \, $, defined in Eq.~(\ref{lambda}), and the relaxation coefficient 
$ \, s_e \,$ used in the relaxation step Eq.~(\ref{moments-collision}), are constant as the spatial step
$ \, \Delta x \,$ tends to zero.
Moreover, we suppose that the relation between the given diffusion coefficient 
$ \, \mu \,$ and the relaxation coefficient $ \, s_J \,$ follows the relation in Eq.~(\ref{mu}). 
In other words, the  relaxation coefficient $ \, s_J \,$ admits the asymptotic
hypotheses in Eqs.~(\ref{s-dx},\ref{developpement-s}) with  $\, s_0 = 0 \,$ 
and $\, s_1 =  {{4+\alpha}\over{6 \, \mu }} $. 
Then, when $ \, \Delta x \,$ tends to zero, the density $ \, \rho \,$ and the  momentum 
after relaxation $ \, J^* \,$ obey the following acoustic model: 
\moneq \label{acoustic-model}  
{{\partial \rho}\over{\partial t}}   \, + \, 
{{\partial J^*}\over{\partial x}}  \,=\, {\rm O}(\Delta x) \,, \quad 
{{\partial J^*}\over{\partial t}}   \, + \, \lambda^2 \, {{4+\alpha}\over{6}} \,   
{{\partial \rho}\over{\partial x}}   \, + \, \lambda^2 \, {{4+\alpha}\over{6 \, \mu}} \, J^*(x, \, t) 
\,=\, {\rm O}(\Delta x) \, . 
\monend

\monitem
The proof of this result is given in  the Appendix~1.
The system Eq.~(\ref{acoustic-model}) is an acoustic model with sound velocity 
 $ \, c_0 = \lambda  \, \sqrt{{4+\alpha}\over{6}} $. 
We see that we have dissipation of momentum with a zero-order operator. 
We have implemented a staggered finite-difference method named
``HaWAY,'' in reference to the authors Harlow and Welch \cite{HaW65}, 
Arakawa \cite{Ar66}    and Yee \cite{Ye65}  who invented it in  the mid 1960's, 
for applications to fluid flow
(``marker and cell''), geophysical sciences (``c-grid'') and electromagnetism 
(``finite difference time domain''), respectively. The details of this second-order 
numerical scheme are given in  Appendix~2. 
This finite-difference approximation gives a correct second-order 
accurate solution of the system obtained by 
replacing the corrections $ \,  {\rm O}(\Delta x) \, $ by $ \, 0 \,$ 
on the right-hand side of Eqs.~(\ref{acoustic-model}).

\bigskip \bigskip   \noindent {\bf \large    5) \quad  
Additional numerical experiments }   

\noindent
We next experiment with the diffusion of a Gaussian density profile with the D1Q3 lattice Boltzmann model
defined in this work. The initial  density profile  is given  by the relation 
\moneq \label{gaussienne-initial} 
 \rho_0(x) =  {\rm exp } \Big( { \displaystyle -{{x^2}\over{4 \, \mu }}} \, \Big) \quad
{\rm with} \quad   \, x \in \R \,. 
\monend
The  other moments $\, J \,$ and $ \, e \,$ are taken to be at equilibrium at $ \, t = 0 $. 
Then the exact solution of the diffusion equation, Eq.~(\ref{edp-diff}), is obtained without difficulty: 
\moneq \label{gaussienne-exact} 
 \rho(x,\, t) =  {{1}\over{\sqrt{1+t}}} \, \exp \Big( -{{x^2}\over{4 \, \mu \, (1+t)}} \,  \Big) 
\,, \quad  x \in \R \,, \ t > 0 \, . 
\monend  
We simulate this problem for $ \, \mu = 0.01 \,$ and $ \, 0 \leq t \leq T=5 \,$ 
in a relatively large domain $ \, -16 \leq x \leq 16  \,$ in order to avoid 
unwanted interactions of the diffusing Gaussian with the boundary.
This has allowed us to employ an elementary periodic boundary condition  
at $  x = \pm 16 $, where all the fields have a value 
inferior to  the smallest number that can be represented in floating-point arithmetic.

\begin{figure}  [H]     \centering
\includegraphics [width=.60 \textwidth, angle=-90] {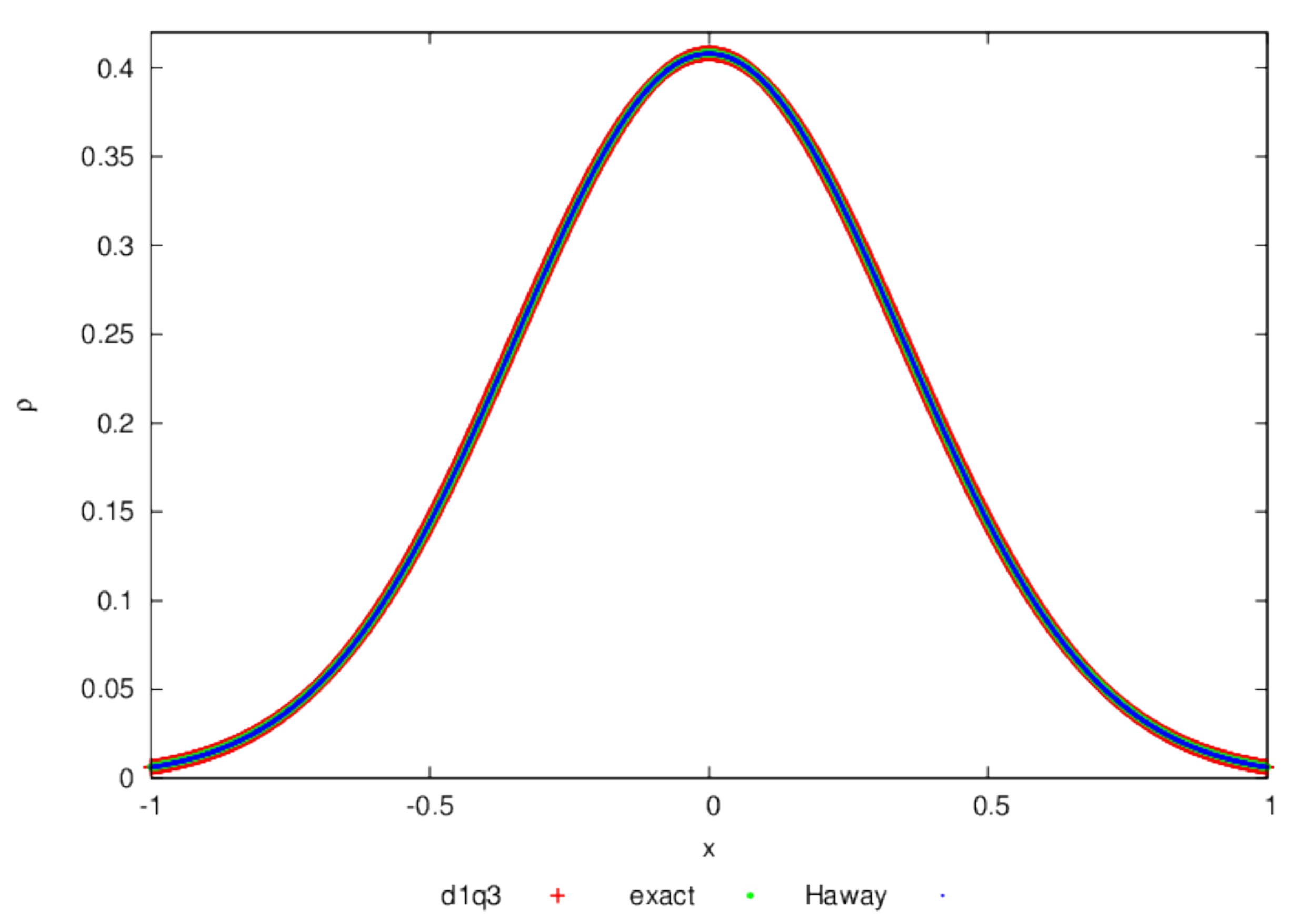}
\caption{Gaussian a time $ \, T = 5  $.}
  \label{gaussian-time5} \end{figure}

\begin{figure}  [H]     \centering
\includegraphics [width=.59 \textwidth, angle=-90] {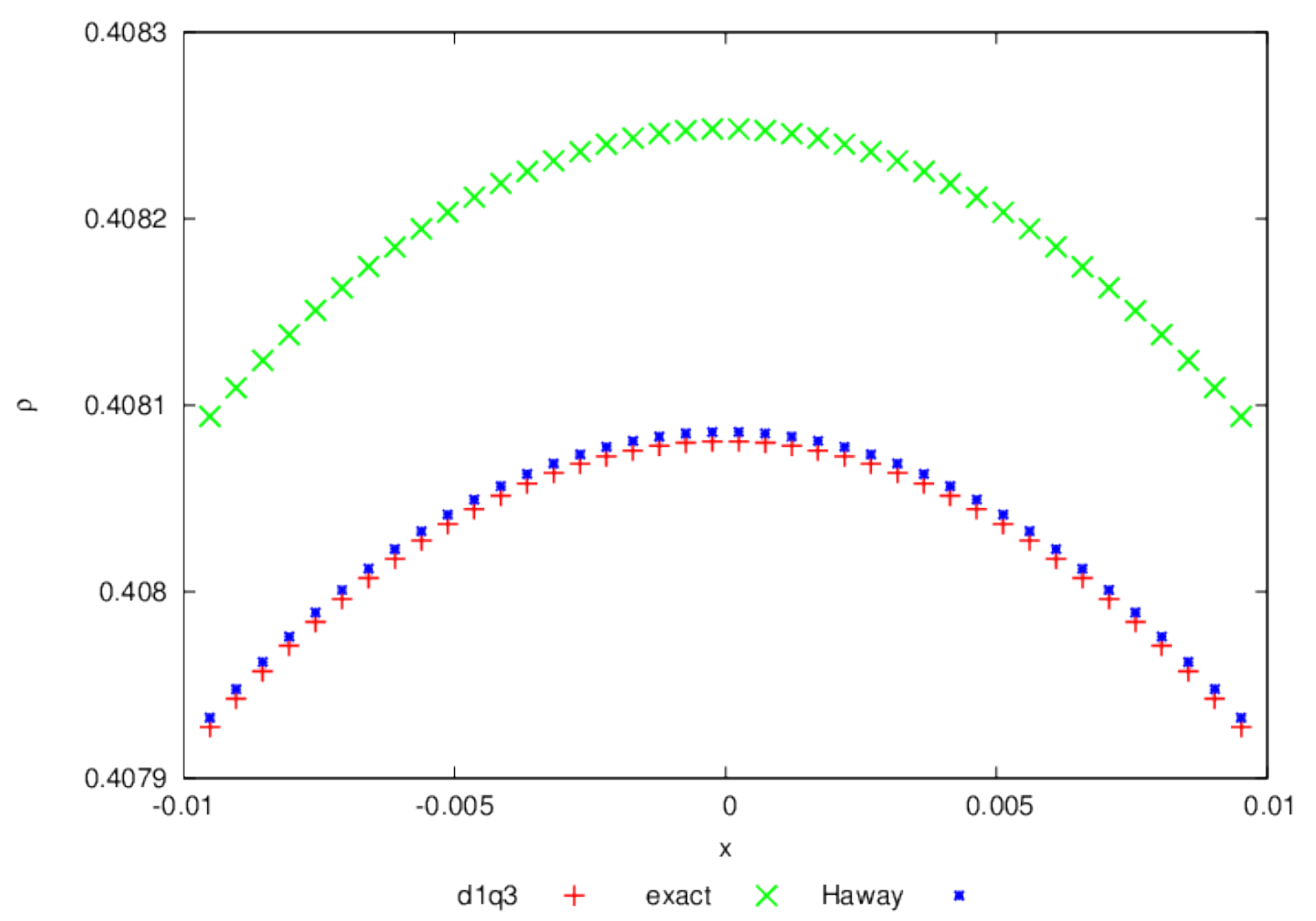}
\caption{Detail of the Gaussian a time $ \, T = 5  $.}
  \label{gaussian-time-zoom2} \end{figure}

\begin{figure}  [H]     \centering
\includegraphics [width=.59 \textwidth, angle=-90] {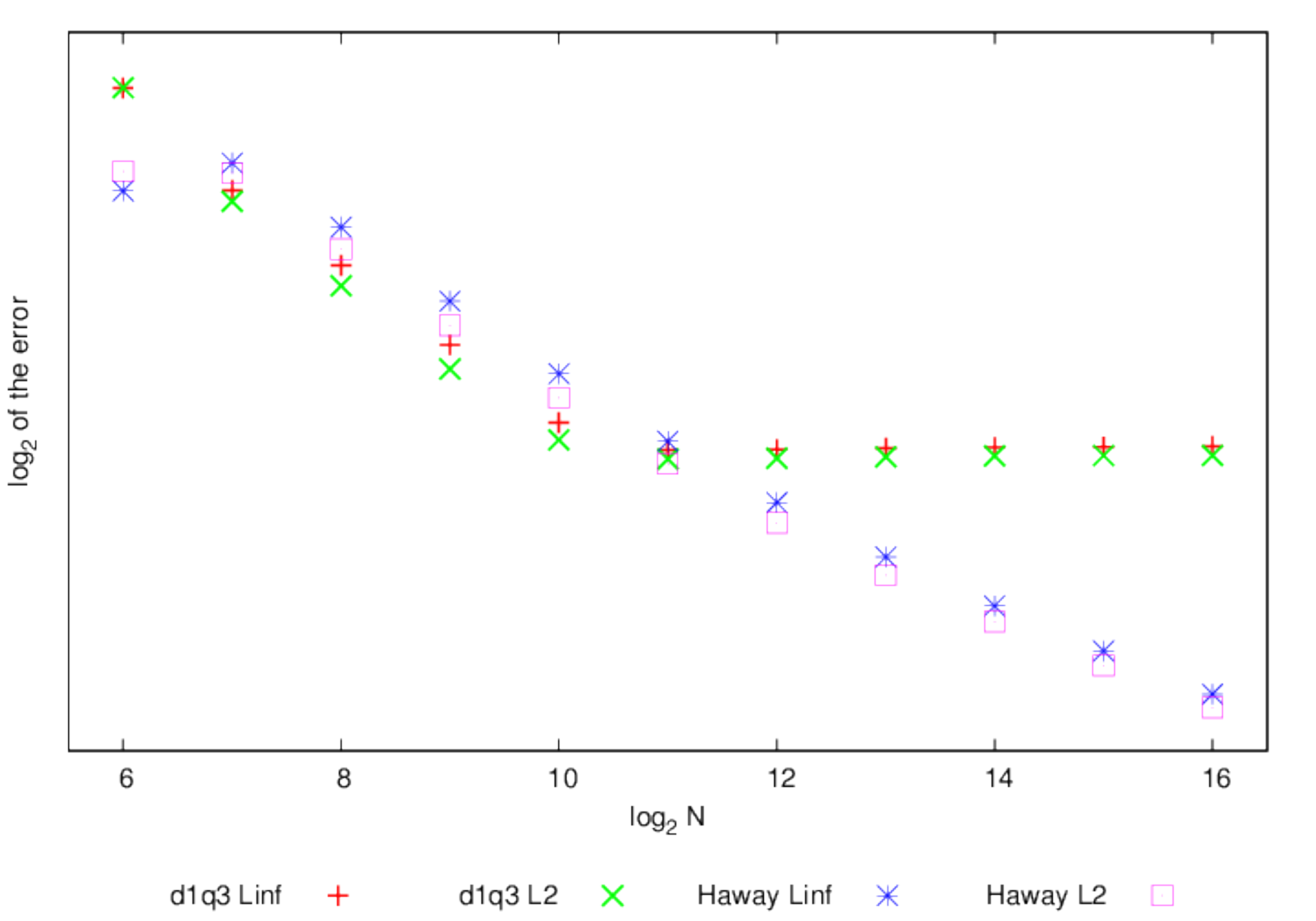}
\caption{D1Q3 lattice Bolztmann scheme for diffusion; Gaussian at time = 5. Non-
convergence towards the exact solution of the diffusion model Eq.~(\ref{edp-diff}) 
and convergence towards the acoustic model Eq.~(\ref{acoustic-model}).} 
  \label{gaussian-non-convergence} \end{figure}

\monitem
At the macroscopic scale, we see in Fig.~\ref{gaussian-time5} that the numerical solution 
computed with the D1Q3 lattice Boltzmann scheme 
faithfully reproduces the exact solution Eq.~(\ref{gaussienne-exact}) of the diffusion equation. 
After magnification by a factor of 100 (Fig.~\ref{gaussian-time-zoom2}), the  D1Q3 model simulates
the acoustic-like system, Eq.~(\ref{acoustic-model}), with better accuracy than it does the diffusion
equation, Eq.~(\ref{edp-diff}).
Fig.~\ref{gaussian-non-convergence} shows that when the mesh is refined 
from $ \, 2^6 = 64 \,$ to $ \, 2^{16} =   65 \, 536 \, $  vertices, the convergence towards the acoustic
model seems reasonable, with an order of accuracy  close to $1$.
In other terms, the diference between the discrete solution of the D1Q3 model 
and the finite-difference simulation of the acoustic model goes  to zero 
 proportionally  to the mesh size.

\begin{figure}  [H]     \centering
\includegraphics [width=.60 \textwidth, angle=-90] {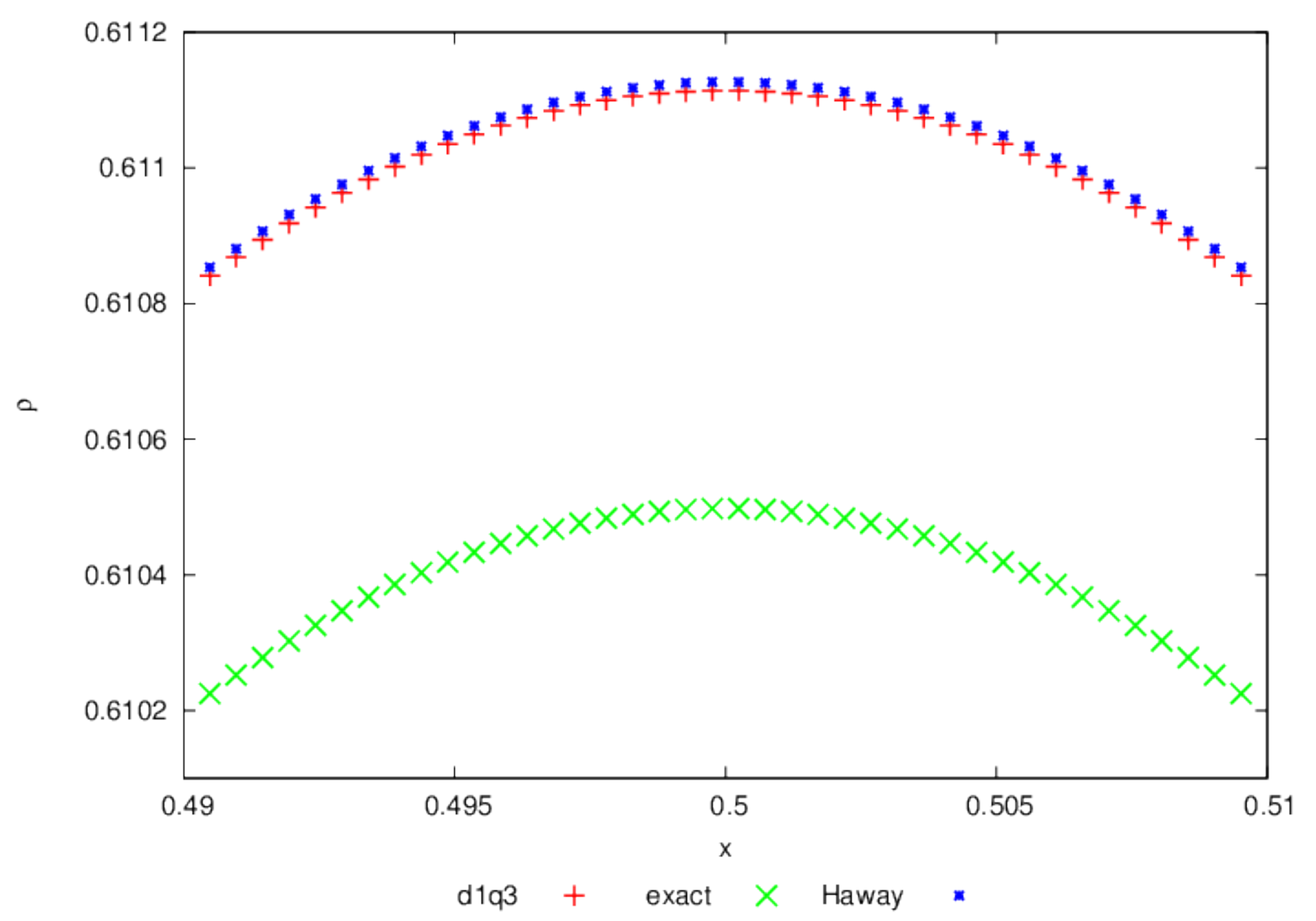}
\caption{Sine wave  at time = 5.  Magnification of the solution around the extremal value. }
\label{sine-time5-second-2}
\end{figure}

\monitem
Now we have what seems like a contradiction:  Our first experiments for the sine wave show
(see, {\it e.g.}, Fig.~\ref{sine-time5-converge})  that 
the diffusion equation is a good reference mathematical model, whereas the acoustic model 
Eq.~(\ref{acoustic-model}) is asymptotically correct for the Gaussian initial condition
(see Fig.~\ref{gaussian-non-convergence}). 
We have performed simulations for the sine wave with much more refined meshes, and lattice sizes
up to 4096. 
At the macroscopic scale, no difference is visible 
between the sine wave solution of the diffusion equation and the 
numerical result proposed by the lattice Boltzmann method 
(again, see Fig.~\ref{sine-time5}). 
After magnification shown in Fig.~\ref{sine-time5-second-2}, the difference between 
the exact solution of the diffusion equation and the D1Q3 solution is 
more important than  the small discrepancy 
between the ``HaWAY'' numerical 
solution of the acoustic model, Eq.~(\ref{acoustic-model}), and the lattice Boltzmann
model. 

\begin{figure}  [H]     \centering
\includegraphics [width=.60 \textwidth, angle=-90] {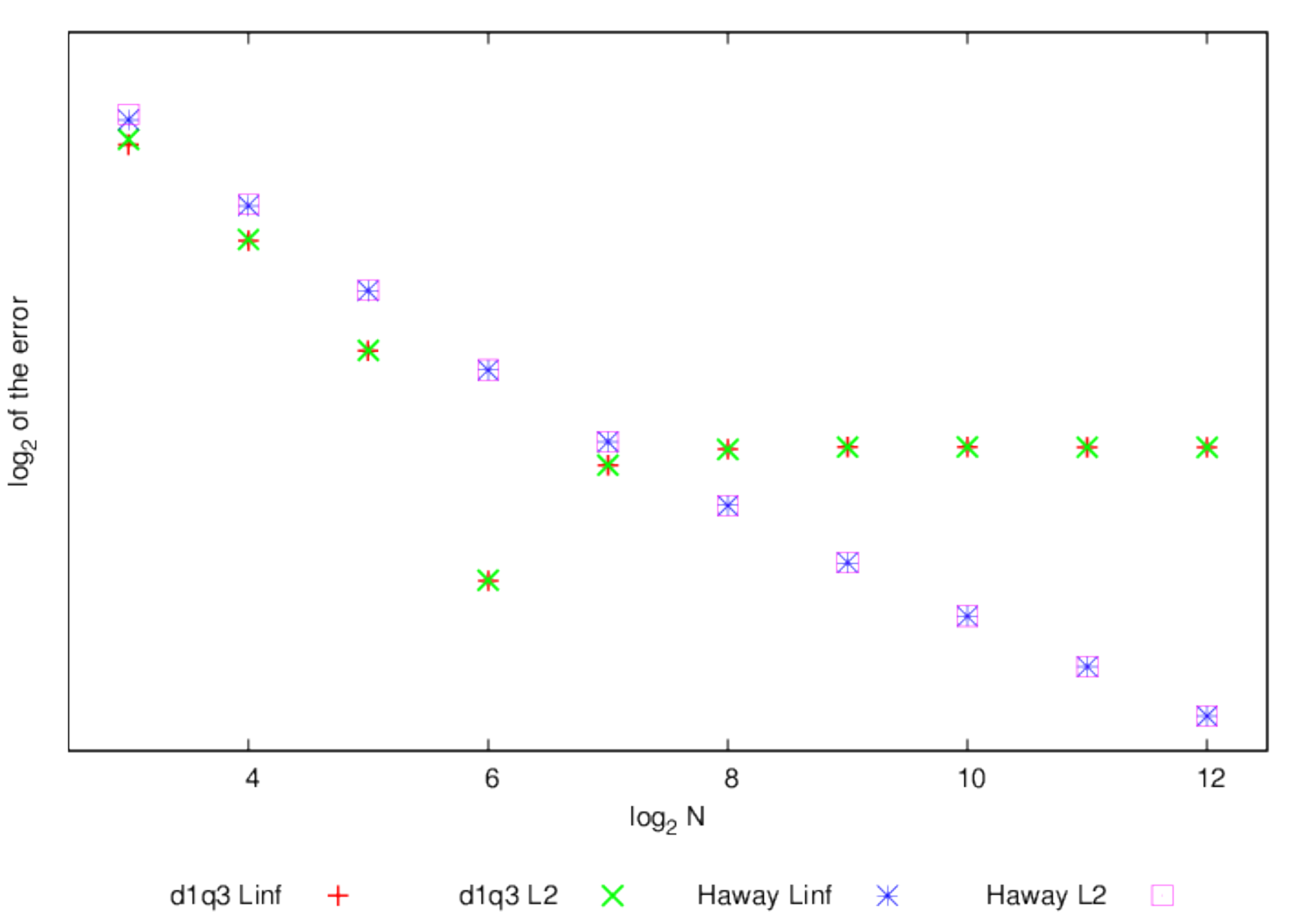}
\caption{D1Q3 lattice Bolztmann scheme for diffusion; Sine wave at time = 5. Non-
convergence towards the exact solution of the diffusion model Eq.~(\ref{edp-diff}) 
and convergence towards the acoustic model Eq.~(\ref{acoustic-model}).}
  \label{sine-non-cv} \end{figure}

\monitem
In Fig.~\ref{sine-non-cv}, we have plotted the quadratic and uniform errors
between the numerical solution obtained from the lattice Boltzmann model 
and the exact solution of the diffusion equation on one hand, and of the approximate solution
(with a second-order scheme) of the acoustic model obtained after a first-order 
Taylor expansion analysis presented at Proposition~1 on the other hand. 
The lattice Boltzmann method gives an excellent approximation of the heat equation 
with the coarse meshes, as shown in Fig.~\ref{sine-time5-converge} in Section~2. 
This good convergence quality cannot be explained by an asymptotic analysis. When the spatial step 
tends to zero, the lattice Boltzmann scheme gives a correct approximation of 
the acoustic model.  Fig.~\ref{sine-non-cv} demonstrates that
the convergence is first-order accurate in both norms.

\bigskip \bigskip   \noindent {\bf \large    6) \quad   Wave vectors analysis }   

\noindent
We may also adopt the point of view of a spectral analysis of the lattice Boltzmann
model, Eqs.~(\ref{moments},\ref{matrice-M},\ref{equilibre},\ref{moments-collision},\ref{schema-temps}). 
We search for a solution of the type
\moneq \label{f-spectral} 
f(x, \, t)  \, = \, \exp \big( i \,  k \, x  \big)  \, \,  \exp \big( \zeta \, t \big)  \, \begin{pmatrix}
\varphi_+ \cr \varphi_0 \cr \varphi_- \end{pmatrix}  \, . 
\monend  
In one time step, we first transform the vector $ \, f \, $ into moments.
Then we relax the moments and return to  the space of particles. 
Finally we advect the result according to Eq.~(\ref{schema-temps}) to recover the particle distribution. 
The collision step $ \, m \longrightarrow m^* \,$ can be written in matrix form: 
\moneqstar 
m^* \, = \, R  \,\, m \,, \quad R \,=\, \begin{pmatrix} 1 & 0 & 0 \cr 
0 & 1-s_J & 0 \cr \alpha \, {{\lambda^2}\over{2}} \, s_e  & 0 & 1 - s_e \end{pmatrix}  
\monendstar 
and the final advection step can be represented by the action of a diagonal matrix:

\moneqstar 
f(x, t + \Delta t)  \, = \, A \, f^* (x, \, t)  \,, \quad
A \,=\, \begin{pmatrix}  \exp ( -i \, \xi )  & 0 & 0 \cr 0 & 1 & 0 \cr 
0 & 0 & \exp ( i \, \xi )  \end{pmatrix}  
\monendstar 
with $ \, \xi =  k \, \Delta x $ . 
Then the vector $ \, \varphi \,$ introduced in Eq.~(\ref{f-spectral}) must be a nontrivial 
solution of the following spectral  problem: 
   \moneqstar 
   \exp \big( \zeta \big) \, \varphi \,=\, 
   A \,\, M^{-1} \,\, R \,\, M \, \varphi  \, , 
   \monendstar   
with $ \, \zeta = z \, \Delta t > 0 $. 
Then, denoting the identity matrix by $ \, {\rm I} \,$, the dispersion relation takes the form 
\moneq \label{dispersion} 
{\rm det}  \, \big[ \,  A \,\, M^{-1} \,\, R \,\, M 
\,-\, \exp \big( \zeta \big) \,  {\rm I} \, \big] \,=\, 0 \,. 
\monend  

\monitem
We have performed  an asymptotic analysis of the relation in Eq.~(\ref{dispersion})
in the limit of a small relaxation coefficient $ \, s_J $
(as in Eq.~(\ref{s-mu})), a small wave vector  $ \, \xi \,$ 
and with a small amplification factor  $ \, \zeta $: 
\moneq \label{asymptotix} 
s_J \,=\, \varepsilon \, s_1 \,,\quad 
\xi \,=\, \varepsilon \, \kappa \,,\quad 
\zeta \,=\, \varepsilon \, \omega \ , 
\monend  
where $ \,  \varepsilon \, $ is a small parameter that tends to zero. 
After some calculation, we obtain without difficulty
\moneq \label{dispersion-eps} 
{\rm det}  \, \big[ \,  A \,\, M^{-1} \,\, R \,\, M \,-\, \exp \big( \zeta \big) \,  {\rm I} \, \big] \,=\, 
- s_e \, \Big( \omega^2  \,+\, s_1 \, \omega \,+\, {{4+\alpha}\over6} \, \kappa \Big) \, \varepsilon^2  
\,+\, {\rm O} (\varepsilon^3) \, . 
\monend  
When $ \, \kappa = 0 $, we recover the hydrodynamic mode with $ \, \omega = 0 \,$
and a dissipative mode according to $ \, \omega = - s_1 $. 
When $ \, \kappa \not= 0 $, we have to solve an equation of degree 2 
made explicit in Eq.~(\ref{dispersion-eps}) at this order of accuracy. 
The discriminant of this equation becomes negative when 
\moneq \label{delta-negatif} 
\kappa \, \geq \, {{s_1}\over{2 \, \sqrt{{4+\alpha}\over6} }}  \, . 
\monend  
%

\begin{figure}  [H]     \centering
\includegraphics [width=.60 \textwidth, angle=-90] {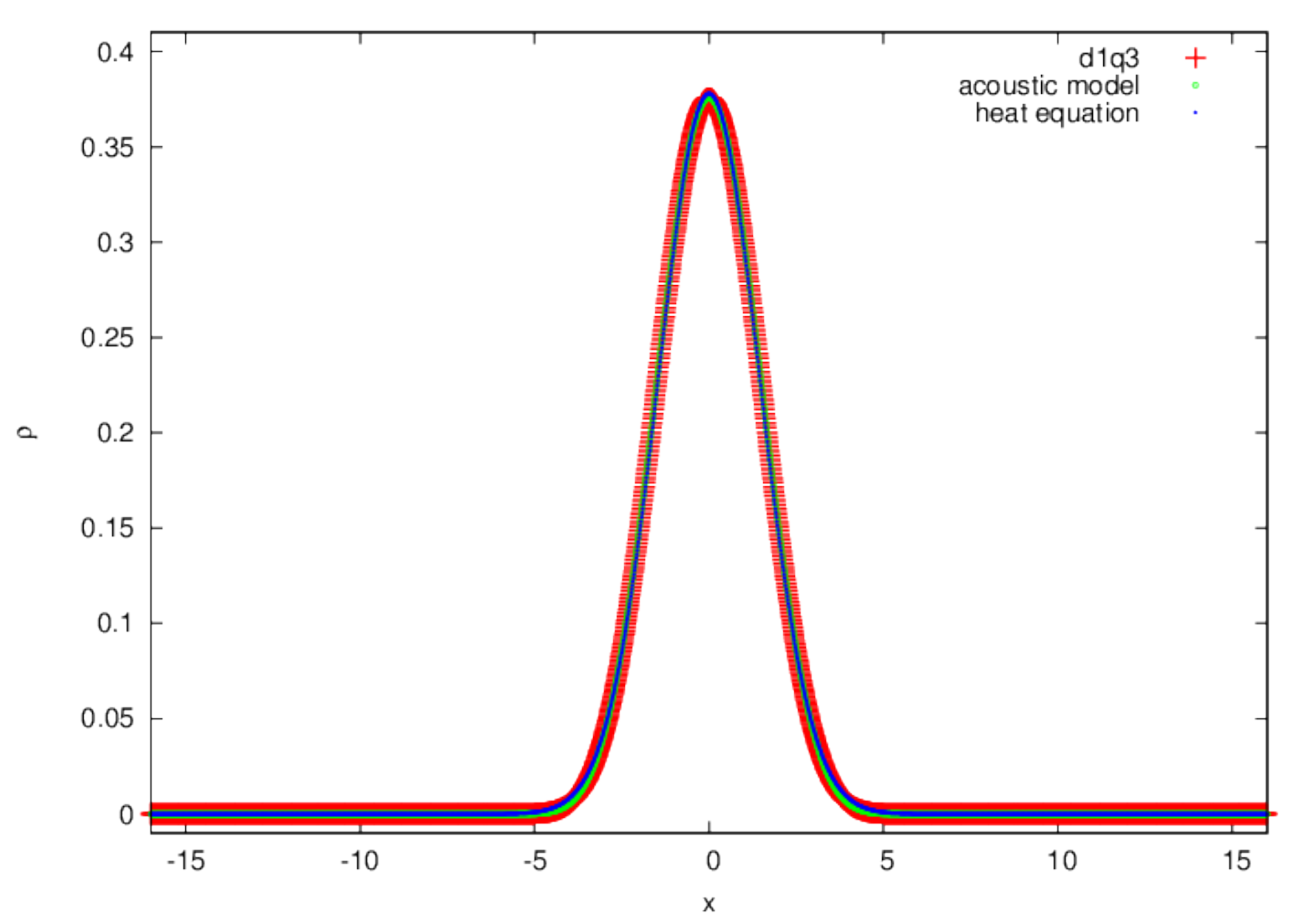}
\caption{Gaussian,  $ \, \mu = 0.15  $. The heat equation is correctly approximated. }
  \label{dromadaire} \end{figure}

\begin{figure}  [H]     \centering
\includegraphics [width=.60 \textwidth, angle=-90] {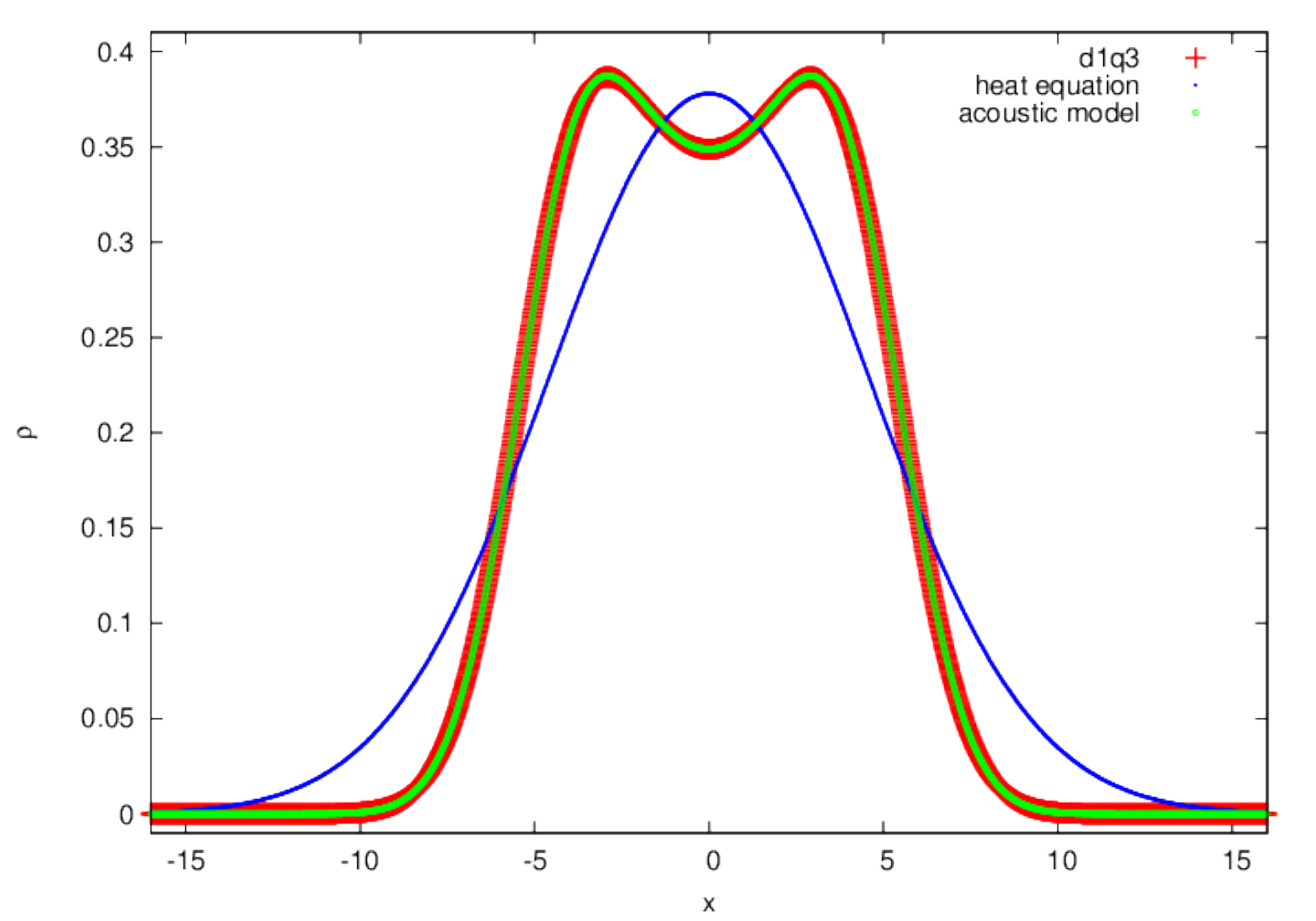}
\caption{Gaussian,  $ \, \mu = 1.5  $. The propagation effects are clearly visible. }
  \label{chameau} \end{figure}

\monitem
We observe that the asymptotics associated with the limit $ \, s_J  \longrightarrow  0 \, $ 
is questionable from a physical point of view. 
When establishing macroscopic partial differential equations it is assumed that internal 
degrees of freedom of the system under study evolve very ``fast'' 
compared to the macroscopic quantities. 
It is known (see, {\it e.g.}, \cite{Du08})
that $ \,  s_J \, $ is given by a ratio of the type 
$ \, {{\Delta t}\over{\tau}}$. 
In the present case, the slow internal degrees of freedom evolve within times 
$ \, \tau \approx {{\Delta t}\over{s}}$. So it is to be expected that the pure diffusion 
partial differential equation will not be accurate for very small values of $s$.

\monitem
In order to make this qualitative difference explicit, we have done two 
numerical experiments with the initial Gaussian profile given by Eq.~(\ref{gaussienne-initial}).
We use $ 2\,048 $  mesh points and a final time  $ \, T=6 \, $ after  384 iterations of the D1Q3 scheme. 
In the first experiment (Fig.~\ref{dromadaire}), 
the diffusion $ \, \mu \,$ is equal to $ \, 0.15 \,$ 
so $ \, s_J = 0.050761 \, $ wich satisfies  Eq.~(\ref{mu}). 
For the second experiment, all parameters are unchanged, except that  $ \, \mu = 1.5 \, $
and the relaxation coefficient  $ \, s_J = 0.005195 $ is much smaller. 
With this value of $ \, s_J \, $  the expression, Eq.~(\ref{delta-negatif}), is satisfied and the discriminant 
of the equation (\ref{dispersion-eps}) is negative.
The propagation effects are evident in Fig.~\ref{chameau}.
One may interpret the result: for small values of parameter $\, s $, the time
between collisions is longer than the duration of the simulation, so particles
move ``ballistically''.

\bigskip  \bigskip    \noindent {\bf \large    7) \quad  
Conclusion }   

We have described a curious convergence property of the D1Q3 lattice Boltzmann model, observed
when trying to simulate a diffusion process with an acoustic scale. 
A new asymptotic analysis has been derived for this circumstance, and we have presented evidence 
of an asymptotic partial differential equation of acoustic type. 
We have observed  analogous difficulties  in two spatial dimensions, 
both for diffusion and Stokes flows. Overall results and physical interpretations
will be given later, with comparison made to the phenomenon of viscoelasticity~\cite{GHL97}. 

A natural question for future study is the generalization of this acoustic-type model to two 
or three spatial dimensions.  Another is the application of this methodology to
lattice Boltzmann models of fluid flow, using an acoustic scale 
while holding fixed the value of the viscosity.

Finally, it seems plausible that there is a link between the strange ``first convergence'' property
noted in this work and the well known tendency of certain asymptotic series to converge at first, followed by divergence~(see {\it e.g.} \cite{Ha49}).  This would raise the question of exactly when the error is minimized, and what is an acceptable approximation of its value when it is minimal.

\bigskip \bigskip   \noindent {\bf  \large  Acknowledgments }   

\noindent 
The authors thank Luc Mieussens who drew our attention to the 
interesting results obtained by St\'ephane Dellacherie in~\cite{De14}.

\bigskip \bigskip   \noindent {\bf \large    Appendix~1  \quad  
Proof of Proposition~1 }    

\noindent
We start from the time iteration, Eq.~(\ref{schema-temps}), and transfer it to the moments:
\moneqstar  
m_k (x ,\, t + \Delta t ) \,=\, 
\sum_{j \, \ell} \, M_{k j} \, \,  M^{-1}_{j \ell}  \,\, m_\ell^* ( x - v_j \, \Delta t , \, t)  \, . 
\monendstar 
With the help of the tensor of momentum-velocity $ \, \Lambda \, $ introduced in~\cite{Du08}, 
defined according to 
\moneqstar 
\Lambda_{k} ^{\ell} \, = \, \sum_{j} \,  M_{k j} \, \, v_j \,\,  M^{-1}_{j \ell} 
\monendstar 
and made explicit for our model as
\moneq \label{tenseur-Lambda} 
\Lambda \, = \, \begin{pmatrix} 0 & 1 & 0 \cr 
{2\over3} \, \lambda^2 & 0 & {1\over3} \cr 0 & \lambda^2 & 0 
\end{pmatrix} \,,  \monend
we have 
\moneq \label{iter-moments} 
m_k (x ,\, t ) \,+\, \Delta t \,  {{\partial m_k}\over{\partial t}}
\,+\, {\rm O}(\Delta t ^2)  \,=\, m_k^* (x ,\, t ) \,-\,  \Delta x \,  \sum_{\ell} \, 
\Lambda_{k} ^{\ell} \,  {{\partial m_\ell^*}\over{\partial x }}  \,+\, {\rm O}(\Delta x ^2) \, . 
\monend

\monitem 
The first moment $ \, m_0 \equiv \rho \,$ is conserved (see Eq.~(\ref{moments-collision})) 
and $ \, \rho^* = \rho $. We deduce from Eq.~(\ref{iter-moments}) and the specific values of
the first line of the matrix $ \, \Lambda \,$ in Eq.~(\ref{tenseur-Lambda}) that
\moneq \label{equation-masse} 
 {{\partial \rho}\over{\partial t}}
\,+\,  {{\partial J^*}\over{\partial x }}  \,=\, {\rm O}(\Delta x) \, 
\monend
and the first equation of Eq.~(\ref{acoustic-model}) is established. 

\monitem 
The third moment $ \, e \, $ is not at equilibrium and we have from the third relation of 
Eqs.~(\ref{moments-collision},\ref{tenseur-Lambda},\ref{iter-moments}):
\moneqstar 
s_e \, \big( e - e^{\rm eq} \big) \, \equiv \, e - e^* \,=\, - \Delta t \,  
{{\partial e}\over{\partial t}} \,-\, 
\Delta x \, \lambda^2 \, {{\partial J^*}\over{\partial x}} 
\,+\,  {\rm O}(\Delta x^2) \,=\,  {\rm O}(\Delta x) \, . 
\monendstar 
The coefficient $ \, s_e \,$ remains constant by hypothesis as $ \, \Delta x \,$ tends to zero. 
Then this moment is close to the equilibrium: 
\moneqstar 
e  \,=\, {{\alpha}\over{2}} \, \lambda^2 \, \rho  \,+\,   {\rm O}(\Delta x) \,,\quad 
e^*  \,=\, {{\alpha}\over{2}} \, \lambda^2 \, \rho  \,+\,   {\rm O}(\Delta x) \,  
\monendstar 
and 
\moneq   \label{moment-e-star} 
{{\partial e^*}\over{\partial x}} \,=\,   \lambda^2 \,  {{\alpha}\over{2}} \,
{{\partial \rho}\over{\partial x}} \,+\,   {\rm O}(\Delta x) \, . 
\monend

\monitem 
The analysis for the second equation differs from what has been proposed previously in~\cite{Du08}
because the moment $ \, J \,$ and the same moment  $ \, J^* \,$ after relaxation are now not close
to the equilibrium value $ \, J^{\rm eq} = 0$. More precisely, we have, due to the second relation of 
Eq.~(\ref{moments-collision}): 
\moneqstar 
J \,=\, {{J^*}\over{1 - s_J}} \,=\,  \Big( 1 \,+\,  {{4+\alpha}\over{6 \, \mu }} \,  \lambda \,  \Delta x 
\,+\, {\rm O}(\Delta x^2) \Big) \, J^*  \, . 
\monendstar 
Then $ \displaystyle \quad 
J(x, \, t + \Delta t) \,=\, \Big( 1 \,+\,  {{4+\alpha}\over{6 \, \mu }} \,  \lambda \,  \Delta x 
\,+\, {\rm O}(\Delta x^2) \Big) \, J^* (x, \, t + \Delta t) $

\smallskip \noindent \qquad  \qquad \qquad  \qquad \quad $ \ \displaystyle 
 \,=\, \Big( 1 \,+\,  {{4+\alpha}\over{6 \, \mu }} \,  \lambda \,  \Delta x  \Big) \, J^* 
 \,+\,  \Delta t \,  {{\partial J^*}\over{\partial t}} \,+\, {\rm O}(\Delta x^2) $ 

\smallskip \noindent \qquad  \qquad \qquad  \qquad \quad $ \ \displaystyle 
 \,=\,  J^* (x, \, t) \,+\,  \Delta t \,  {{\partial J^*}\over{\partial t}}  
\,+\,  {{4+\alpha}\over{6 \, \mu }} \,  \lambda \,  \Delta x  \, J^*(x, \, t) 
  \,+\, {\rm O}(\Delta x^2)  \, . $ 

\smallskip \noindent 
We report this expression in the expansion Eq.~(\ref{iter-moments}), we substract 
$ \, J^*(x, \, t) \,$ from both sides of the equation and we divide by $ \, \Delta t$. 
Due to the previous result Eq.~(\ref{moment-e-star}), we obtain:  
\moneqstar 
 {{\partial J^*}\over{\partial t}}  
\,+\,  {{4+\alpha}\over{6 \, \mu }} \,  \lambda^2 \,  J^*  \,+\, {\rm O}(\Delta x) 
\,=\, - {2\over3} \, \lambda ^2 \, {{\partial \rho}\over{\partial x}} \, - \, 
\lambda^2 \,  {{\alpha}\over{2}} \, {{\partial \rho}\over{\partial x}} \,+\,   {\rm O}(\Delta x) 
\monendstar 
and the second equation of Eqs.~(\ref{acoustic-model})  is established. \hfill $\square$

\bigskip \bigskip     \noindent {\bf \large    Appendix~2  \quad  
``HaWAY'' staggered finite differences }    

\begin{figure}  [H]     \centering
\includegraphics [width=.51 \textwidth, angle=0] {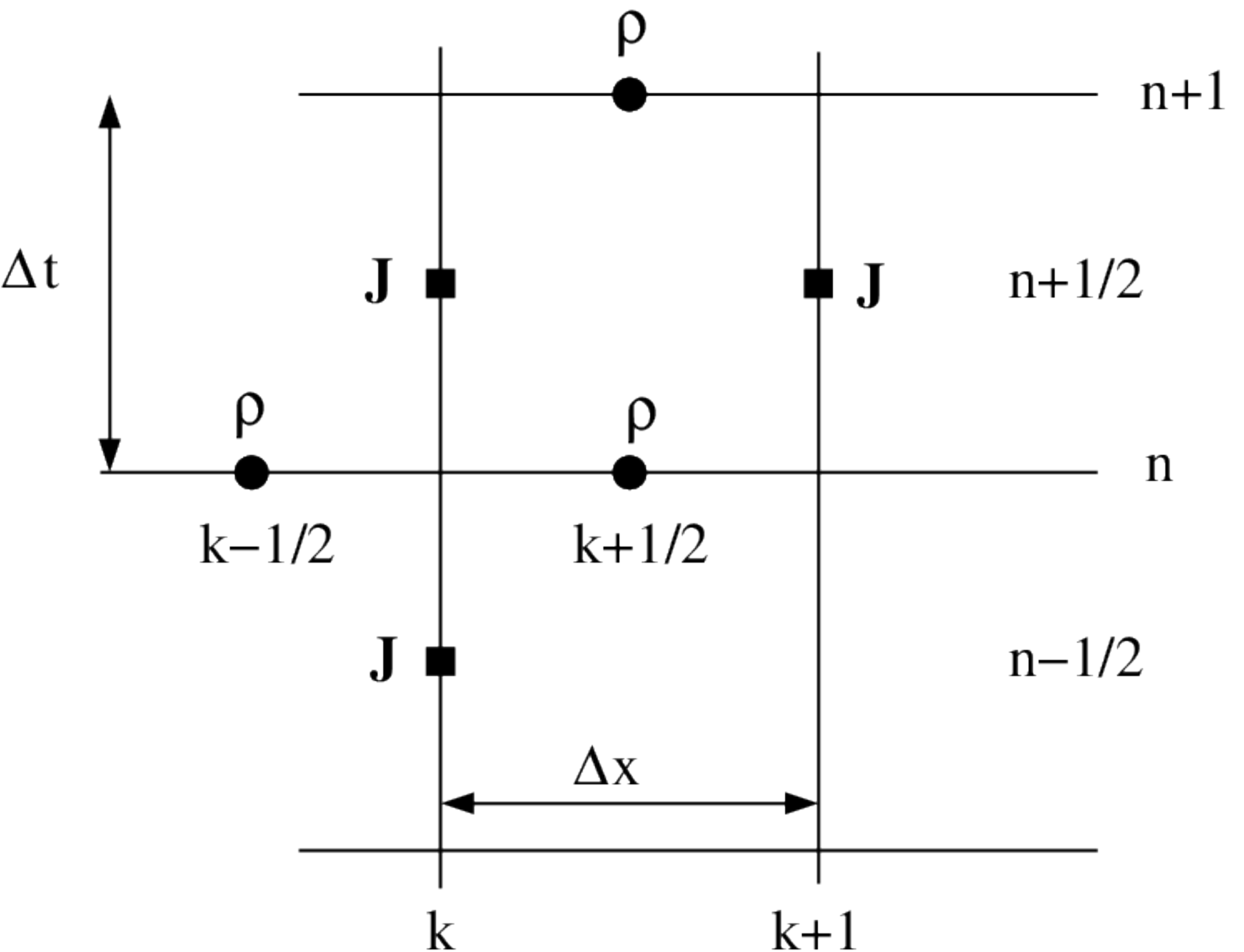}
\caption{HaWAY grid for staggered finite differences.  }
  \label{grille-HaWAY} \end{figure}

\noindent 
We consider the acoustic model proposed in Eq.~(\ref{acoustic-model}).  With compact notation, 
we denote it here according to~:
\moneq \label{acoustic-model-compact}  
{{\partial \rho}\over{\partial t}}   \, + \, 
{{\partial J}\over{\partial x}}  \,=\, 0  \,, \quad 
{{\partial J}\over{\partial t}}   \, + \, c_0^2 \, {{\partial \rho}\over{\partial x}}   
\, + \, \Gamma \,  J(x, \, t)   \,=\, 0  \, . 
\monend
Given a spatial step $ \, \Delta x \, $  and a time step $ \, \Delta t  $, 
we consider  integer multiples of these parameters for the discretization of space and time. 
The density $ \,  \rho  \,$ is approximated at semi-integer vertices in space and integer points in time
whereas the momentum  $ \,  J  \,$  is approximated at integer nodes  in space 
and semi-integer values in time:
\moneq \label{variables-discretes}  
 \rho \, \approx \, \rho ^n_{k+1/2} \, , \quad J  \, \approx \, J^{n+1/2}_k \, . 
\monend
The Figure~\ref{grille-HaWAY} gives an illustration of this classical choice 
\cite{Ar66, HaW65, Ye65}. 

\monitem 
We discretize the first equation of Eqs.~(\ref{acoustic-model-compact}) with a two-point centered 
finite-difference schemes around the vertex 
$ \, \big( (k+{1\over2}) \, \Delta x \,,\, (n+{1\over2}) \, \Delta t \big) $: 
\moneq \label{masse-discrete}  
 {{1}\over{\Delta t}} \, \big( \rho^{n+1}_{k+1/2} - \rho ^n_{k+1/2} \big) \,+\, 
 {{1}\over{\Delta x}} \, \big( J^{n+1/2}_{k+1} - J^{n+1/2}_{k} \big)  \,=\, 0  \, . 
\monend
We use the same approach for the discretization of the second  equation of Eqs.~(\ref{acoustic-model-compact}) 
 around the node  
$ \, \big( k \, \Delta x \,,\, n \, \Delta t \big) $:  
\begin{equation} \label{impulsion-discrete-1}  \displaystyle 
 {{1}\over{\Delta t}} \, \big(  J^{n+1/2}_{k} -  J^{n-1/2}_{k} \big) \,+\, 
 {{1}\over{\Delta x}} \, \big( \rho^{n}_{k+1/2}  - \rho^{n}_{k-1/2}  \big) 
\,+\, \Gamma \, J^{n}_{k} \,=\, 0  \, . 
\end{equation}
We interpolate the momentum at integer vertices  with a simple centered mean value: 
\moneqstar 
  J^{n}_{k} \,=\, {{1}\over{2}} \, \big(  J^{n+1/2}_{k} +  J^{n-1/2}_{k} \big)  \, . 
\monendstar 
We incorporate this expression into  the relation Eq.~(\ref{impulsion-discrete-1}) 
and we obtain 
\moneq \label{impulsion-discrete-2}  
\Big(  {{1}\over{\Delta t}} + {{\Gamma}\over{2}} \, \Big) \,   J^{n+1/2}_{k} 
 \,+\,  {{1}\over{\Delta x}} \, \big( \rho^{n}_{k+1/2}  - \rho^{n}_{k-1/2}  \big) 
\,=\, 
\Big(  {{1}\over{\Delta t}} - {{\Gamma}\over{2}} \, \Big) \, J^{n-1/2}_{k} \, . 
\monend
The numerical scheme is now entirely defined for internal nodes. 
We have used periodic boundary conditions.

\bigskip \bigskip      \noindent {\bf  \large  References }   


\end{document}